\documentclass[12pt]{amsart}
\usepackage{latexsym, amssymb, amsmath, amscd}
 \usepackage{eucal}

\theoremstyle{plain}
    \newtheorem{rema}{Remark}[section]

   \newtheorem{theo}[rema]{Theorem}
 \newtheorem{conj}[rema]{Conjecture}
   
    \newtheorem{lemma}[rema]{Lemma}
    \newtheorem{corol}[rema]{Corollary}
     
  \newtheorem{rmk}[rema]{Remark}

	\newcommand{\nno}{\nonumber}

	\newcommand{\p}{\partial}
\newcommand{\fr}{\frac}

 \newcommand{\pf}{{\it Proof:}\hspace{2ex}}
 
 \newcommand{\epfv}{\hspace{1em}$\Box$\vspace{1em}}

\newcommand{\bC}{{\mathbb C}}

\newcommand{\bN}{{\mathbb N}}

\newcommand{\BQ}{\begin{eqnarray}}
\newcommand{\EQ}{\end{eqnarray}}
\newcommand{\BQn}{\begin{eqnarray*}}
\newcommand{\EQn}{\end{eqnarray*}}

\newcommand{\BL}{\begin{align}}
\newcommand{\EL}{\end{align}}
\newcommand{\BLn}{\begin{align*}}
\newcommand{\ELn}{\end{align*}}

\renewcommand{\theequation}{\thesection.\arabic{equation}}
\renewcommand{\therema}{\thesection.\arabic{rema}}
\newcommand{\EAn}{\end{align*}}

\newcommand{\Hes}{ \text{Hes\,} }

\title[Homogeneous Hessian Nilpotent Polynomials]
{Two Results on Homogeneous Hessian Nilpotent Polynomials}

    \author[Arno van den Essen and Wenhua Zhao]
{Arno van den Essen$^{*}$ and Wenhua Zhao$^{**}$}      
    
\begin{document}

\begin{abstract}
Let $z=(z_1,  \cdots, z_n)$ and 
$\Delta=\sum_{i=1}^n \fr {\p^2}{\p z^2_i}$
the Laplace operator. A formal power series 
$P(z)$ is said to be
{\it Hessian Nilpotent}(HN) 
if its Hessian matrix $\Hes P(z)=(\fr {\p^2 P}{\p z_i\p z_j})$ 
is nilpotent. In recent developments 
in \cite{BE1}, \cite{M} and \cite{HNP},
the Jacobian conjecture 
has been reduced to the following so-called  
{\it vanishing conjecture}(VC) 
of HN polynomials: {\it for any homogeneous HN polynomial 
$P(z)$ $($of degree $d=4$$)$, we have
$\Delta^m P^{m+1}(z)=0$ for any 
$m>>0$.} In this paper, 
we first show that, the VC holds 
for any homogeneous HN  
polynomial $P(z)$ provided that 
the projective subvarieties ${\mathcal Z}_P$ and 
${\mathcal Z}_{\sigma_2}$ of $\bC P^{n-1}$ 
determined by the principal ideals 
generated by $P(z)$ and $\sigma_2(z)\!:=\sum_{i=1}^n z_i^2$, 
respectively, intersect only 
at regular points of ${\mathcal Z}_P$.
Consequently, the Jacobian conjecture holds 
for the symmetric polynomial maps 
$F=z-\nabla P$ with $P(z)$ HN if $F$ has 
no non-zero fixed point $w\in \bC^n$ with $\sum_{i=1}^n w_i^2=0$.
Secondly, we show that the VC holds 
for a HN formal power series $P(z)$ if and only if,  
for any polynomial $f(z)$, $\Delta^m (f(z)P(z)^m)=0$ 
when $m>>0$.
\end{abstract}

\keywords{Hessian nilpotent polynomials, 
the vanishing conjecture, symmetric polynomial maps, 
the Jacobian conjecture.}
   
\subjclass[2000]{14R15, 31B05.}

 \bibliographystyle{alpha}
    \maketitle


\renewcommand{\theequation}{\thesection.\arabic{equation}}
\renewcommand{\therema}{\thesection.\arabic{rema}}
\setcounter{equation}{0}
\setcounter{rema}{0}
\setcounter{section}{0}

\section{\bf Introduction and Main Results}

Let $z=(z_1, z_2, \cdots, z_n)$ be 
commutative free variables.
Recall that the well-known 
Jacobian conjecture claims that:
{\it  any polynomial map 
$F(z): \bC^n \to \bC^n$ with the Jacobian 
$j(F)(z)\equiv 1$ is an autompophism of $\bC^n$
and its inverse map must also 
be a polynomial map}. 
Despite intense study 
from mathematicians in more 
than sixty years, 
the conjecture is still 
open even for the case $n=2$. 
In 1998, S. Smale \cite{S} 
included the Jacobian conjecture 
in his list of $18$ important 
mathematical problems 
for $21$st century.
For more history and known results on
the Jacobian conjecture, 
see \cite{BCW}, \cite{E} 
and references there. 

Recently, M. de Bondt and the first author
\cite{BE1} and G. Meng \cite{M} independently 
made the following remarkable breakthrough 
on the Jacobian conjecture. Namely,
they reduced the Jacobian conjecture to the 
so-called {\it symmetric} polynomial maps, 
i.e the polynomial maps of the form
$F=z-\nabla P$, 
where $\nabla P\! :=(\fr{\p P}{\p z_1}, \fr{\p P}{\p z_2}, 
\cdots ,\fr{\p P}{\p z_n})$, i.e. $\nabla P(z)$ 
is the {\it gradient} of $P(z) \in \bC[z]$.

For more recent developments on the Jacobian 
conjecture for symmetric polynomial maps, 
see \cite{BE1}--\cite{BE4}.

Based on the symmetric reduction 
above and also the classical 
homogeneous reduction in \cite{BCW} 
and \cite{Y}, the second author in \cite{HNP} 
further reduced the Jacobian conjecture to 
the following so-called vanishing conjecture.

Let $\Delta\! :=\sum_{i=1}^n \fr {\p^2}{\p z^2_i}$
the Laplace operator and call a formal power series 
$P(z)$ {\it Hessian nilpotent}(HN) if its 
Hessian matrix $\Hes P(z)\! :=(\fr {\p^2 P}{\p z_i\p z_j})$ 
is nilpotent. It has been shown in \cite{HNP} 
that the Jacobian conjecture is equivalent to

\begin{conj} {\bf (Vanishing Conjecture of HN Polynomials)} \\
{\it For any homogeneous HN polynomial $P(z)$ $($of 
degree $d=4$$)$, we have $\Delta^m P^{m+1}=0$ when $m>>0$.}
\end{conj}

Note that, it has also been shown 
in \cite{HNP} that $P(z)$ is HN 
if and only if $\Delta^m P^{m}=0$ 
for $m\geq 1$. 
 
In this paper, we will prove the following 
two results on HN polynomials.

Let $P(z)$ be a homogeneous HN polynomial of degree 
$d\geq 3$ and $\sigma_2(z)\!\! :=\!\sum_{i=1}^n z_i^2$. 
We denote by ${\mathcal Z}_P$ and ${\mathcal Z}_{\sigma_2}$ 
the projective subvarieties of $\bC P^{n-1}$ 
determined by the principal ideals generated by 
$P(z)$ and $\sigma_2(z)$, respectively.
The first main result of this paper 
is the following theorem.

\begin{theo}\label{MainResult-1}
Let $P(z)$ be a homogeneous HN polynomial
of degree $d \geq 4$. Assume that 
${\mathcal Z}_P$ intersects with 
${\mathcal Z}_{\sigma_2}$ only 
at regular points of ${\mathcal Z}_P$, 
then the vanishing conjecture 
holds for $P(z)$. In particular, 
the vanishing conjecture holds 
if the projective variety ${\mathcal Z}_P$ 
is regular. 
\end{theo}

\begin{rmk}
Note that, when $\deg P(z)=d=2$ or $3$, 
the Jacobian conjecture holds 
for the symmetric polynomial map $F=z-\nabla P$. 
This is because, when $d=2$, $F$ is a linear map with 
$j(F)\equiv 1$. Hence $F$ is an automorphism of $\bC^n$; 
while when $d=3$, we have $\deg F=2$. By Wang's theorem 
\cite{Wa}, the Jacobian conjecture holds for $F$ again. 
Then, by the equivalence of the vanishing conjecture 
for the homogeneous HN polynomial $P(z)$ 
and the Jacobian conjecture for the symmetric map 
$F=z-\nabla P$ established in \cite{HNP}, we see that,
when $\deg P(z)=d=2$ or $3$, Theorem \ref{MainResult-1} 
actually also holds even without 
the condition on the projective variety 
${\mathcal Z}_P$. 
\end{rmk}

For any non-zero $z\in \bC^n$, denote 
by $[z]$ its image in the projective 
space $\bC P^{n-1}$. Set 
\begin{align}
\widetilde {\mathcal Z}_{\sigma_2}\! :=\{z \in \bC^n \,|\, z\neq 0;\, 
[z] \in {\mathcal Z}_{\sigma_2}\}. 
\end{align}
In other words, $\widetilde {\mathcal Z}_{\sigma_2}$ is the set 
of non-zero $z\in \bC^n$ such that $\sum_{i=1}^n z_i^2=0$.

Note that, for any homogeneous polynomial 
$P(z)$ of degree $d$, it follows from the Euler's 
formula $dP=\sum_{i=1}^n z_i\frac{d P}{d z_i}$, 
that any non-zero $w\in \bC^n$, 
$[w]\in \bC P^{n-1}$ is a singular point 
of ${\mathcal Z}_P$ if and only if $w$ 
is a fixed point of the symmetric map 
$F=z-\nabla P$. Furthermore, 
it is also well-known that, 
$j(F)\equiv 1$ if and only if $P(z)$ is HN.

By the observations above and Theorem \ref{MainResult-1}, 
it is easy to see that we have the following corollary 
on symmetric polynomial maps.

\begin{corol}\label{MR1-corol-1}
Let $F=z-\nabla P$ with $P$ 
homogeneous and $j(F)\equiv 1$ 
$($or equivalently, $P$ is HN$)$. 
Assume that $F$ does not fix any 
$w\in \widetilde {\mathcal Z}_{\sigma_2}$. 
Then the Jacobian holds for $F(z)$.  
In particular, if $F$ has no 
non-zero fixed point, 
the Jacobian conjecture holds for 
$F$.
\end{corol}

Our second main result is following theorem which says that 
the vanishing conjecture is actually equivalent to a formally 
much stronger statement. 

\begin{theo}\label{MainResult-2}
For any HN polynomial $P(z)$, the vanishing conjecture 
holds for $P(z)$ if and only if, for any polynomial 
$f(z)\in \bC[z]$, $\Delta^m (f(z)P(z)^m)=0$ 
when $m>>0$.
\end{theo}

\section{\bf Proof of the Main Results}

Let us first fix the following notation. 
Let $z=(z_1, z_2,\cdots, z_n)$
be free complex variables
and $\bC[z]$ (resp.\,$\bC[[z]]$) the algebra of 
polynomials (resp.\,formal power series) in $z$.  
For any $d \geq 0$, we denote by $V_d$ the vector 
space of homogeneous polynomials in $z$ of degree $d$.

For any $1\leq i\leq n$, we set 
$D_i=\fr {\p}{\p z_i}$ 
and $D=(D_1, D_2,\cdots, D_n)$.
We define a $\bC$-bilinear map 
$\{\cdot, \cdot\}: \bC[z] \times \bC[z] \to   \bC[z]$ by setting
\begin{align*}
\{f,  \, g \}\! := f(D)g(z)
\end{align*}
for any $f(z), g(z)\in \bC[z]$.

Note that, for any $m \geq 0$, 
the restriction of $\{\cdot, \cdot\}$
on $V_m \times V_m $ gives a
$\bC$-bilinear form of the 
vector subspace $V_m$, 
which we will denote by $B_m(\cdot, \cdot)$.
It is easy to check that, for any $m\geq 1$,
$B_m(\cdot, \cdot)$ is symmetric 
and non-singular.

The following lemma will play 
a crucial role 
in our proof of the first main result. 

\begin{lemma}\label{MainLemma}
For any homogeneous polynomials 
$g_i(z)$ $(1\leq i\leq k)$ of degree $d_i\geq 1$, let $S$ be the vector 
space of polynomial solutions of the following system of PDEs: 
\begin{align}\label{Sys-1}
\begin{cases}
g_1 (D)\, u(z)=0, \\
g_2 (D)\, u(z)=0, \\
\quad ..... \quad \\
g_k (D)\, u(z)=0.
\end{cases}
\end{align}
Then, $\dim S<+\infty$ 
if and only if $g_i(z)$ $(1\leq i\leq k)$ 
have no non-zero common zeroes.
\end{lemma}
 
\pf
Let $I$ the homogeneous ideal of $\bC[z]$ 
generated by $\{ g_i(z) | 1\leq i\leq k\}$. 
Since all $g_i(z)$'s are homogeneous, 
$S$ is a homogeneous 
vector subspace $S$ of $\bC[z]$. 

Write 
\begin{align}
S& =\bigoplus_{m=0}^\infty S_m, \\
I& =\bigoplus_{m=0}^\infty I_m.
\end{align}
where $I_m\! :=I\cap V_m$ and 
$S_m\! :=I\cap V_m$ for any $m\geq 0$. 

\vskip2mm

\underline{\it Claim:} {\it For any 
$m\geq 1$ and $u(z)\in V_m$, 
$u(z)\in S_m$ if and only if $\{u, I_m \}=0$, or 
in other words, $S_m=I_m^\perp$ 
with respect to the $\bC$-bilinear form 
$B_m(\cdot, \cdot)$ of $V_m$.}

\vskip2mm

\underline{\it Proof of the Claim:}
First, by the definitions of $I$ and $S$, 
we have $\{I_m, S_m\}=0$ for any $m\geq 1$, 
hence $S_m\subseteq I_m^\perp$. Therefore, 
we need only show that, 
for any $u(z)\in I_m^\perp\subset V_m$, 
$g_i(D)u(z)=0$ for any $1\leq i\leq n$. 

We first fix any $1\leq i\leq n$. 
If $m < d_i$, there is nothing to prove. 
If $m=d_i$, then $g_i(z) \in I_m$, 
hence $\{g_i, u\}=g_i(D)u=0$.
Now suppose $m>d_i$. Note that, for any 
$v(z)\in V_{m-d_i}$, $v(z) g_i(z)\in I_m$. 
Hence we have 
\begin{align*}
0&= \{ v(z) g_i(z), u(z) \}\\
&=v(D)g_i(D) u(z)\\
&=v(D)\left (g_i(D) u \right )(z)\\
&=\{v(z),  \left (g_i(D)u\right )(z) \}.
\end{align*}

Therefore, we have 
$$
B_{m-d_i}\left( (g_i(D)u) (z),\, V_{m-d_i} \right )=0.
$$
Since $B_{m-d_i}(\cdot, \cdot)$ is a 
non-singular $\bC$-bilinear form of $V_{m-d_i}$, 
we have $g_i(D)u=0$. Hence, the Claim holds. \epfv


By a well-known fact in Algebraic Geometry 
(see Exercise $2.2$ in \cite{H}, for example), 
we know that the homogeneous polynomials
$g_i(z)$ $(1\leq i\leq k)$ have no non-zero 
common zeroes if and only if 
$I_m=V_m$ when $m>>0$. While, by the Claim above, 
we know that, $I_m=V_m$ when $m>>0$ if and only if $S_m=0$ when $m>>0$, 
and if and only if the solution space $S$ of the system (\ref{Sys-1}) 
is finite dimensional. 
Hence, the lemma follows.
\epfv 

Now we are ready to prove our first main result, 
Theorem \ref{MainResult-1}.

\vskip2mm

\underline{\it Proof of Theorem \ref{MainResult-1}}: 
Let $P(z)$ be a homogeneous HN polynomial 
of degree $d\geq 4$ and $S$ the vector space 
of polynomial solutions 
of the following system of PDEs:
\begin{align}\label{Sys-2}
\begin{cases}
\frac{\p P}{\p z_1} (D)\, u(z)=0, \\
\frac{\p P}{\p z_2} (D)\, u(z)=0, \\
\quad ..... \quad \\
\frac{\p P}{\p z_n} (D)\, u(z)=0, \\
\Delta \, u(z)=0.
\end{cases}
\end{align}

First, note that the projective subvariety 
${\mathcal Z}_P$ intersects with 
${\mathcal Z}_{\sigma_2}$ only 
at regular points of ${\mathcal Z}_P$ 
if and only if $\frac{\p P}{\p z_i}(z)$ 
$(1\leq i\leq n)$ and 
$\sigma_2=\sum_{i=1}^n z_i^2$ 
have no non-zero common zeros (agian use Euler's formula). 
Then, by Lemma \ref{MainLemma}, 
we have $\dim S < +\infty$.

On the other hand, by Theorem $6.3$ in \cite{HNP}, 
we know that $\Delta^m P^{m+1}\in S$ 
for any $m\geq 0$.
Note that $\deg \Delta^m P^{m+1}=(d-2)m+d$ 
for any $m\geq 0$. So $\deg \Delta^m P^{m+1}>\deg \Delta^k P^{k+1}$ 
for any $m>k$. Since $\dim S < +\infty$ (from above), 
we have $\Delta^m P^{m+1}=0$ when $m>>0$, 
i.e. the vanishing conjecture holds for $P(z)$. 
\epfv

Next, we give a proof for our second main result,
Theorem \ref{MainResult-2}.

\vskip2mm

\underline{\it Proof of Theorem \ref{MainResult-2}}:
The $(\Leftarrow)$ part follows directly 
by choosing $f(z)$ to be $P(z)$ itself.

To show $(\Rightarrow)$ part, let $d=\deg f(z)$. If $d=0$, 
$f$ is a constant. Then, $\Delta^m (f(z)P(z)^m)=f(z)\Delta^m P^m=0$
for any $m\geq 1$.

So we assume $d\geq 1$. By Theorem $6.2$ in \cite{HNP}, 
we know that, if the vanishing conjecture holds 
for $P(z)$, then, for any fixed $a\geq 1$, 
$\Delta^m P^{m+a}=0$ when $m>>0$. 
Therefore there exists 
$N>0$ such that, for any $0\leq b\leq d$ 
and any $m>N$, we have $\Delta^m P^{m+b}=0$.

By Lemma $6.5$ in \cite{HNP}, 
for any $m\geq 1$, we have
\begin{align}\label{MainResult-2-pe1}
&\Delta^m (f(z) P(z)^m)= \\
&\, \sum_{\substack{k_1+ k_2+ k_3= m\\
k_1, k_2, k_3 \geq 0}} 2^{k_2}\binom {m}{k_1, k_2, k_3}
\sum_{\substack{ {\bf s} \in \bN^n \\ |{\bf s}|=k_2}} 
\binom {k_2}{\bf s} \frac {\p^{k_2} \Delta^{k_1} f(z)}{\p z^{\bf s}}
\frac {\p^{k_2}\Delta^{k_3} P^m(z) }{\p z^{\bf s}},\nno 
\end{align}
where $\binom {m}{k_1, k_2, k_3}$ and $\binom{k_2}{\bf s}$ 
denote the usual binomials.

Note first that, the general term in the sum above is non-zero only 
if $2k_1+k_2 \leq d$. But on the 
other hand, since 
\begin{align}\label{MainResult-2-pe2}
0\leq k_1+k_2\leq 2k_1+k_2\leq d,
\end{align}
by the choice of $N\geq 1$, we have
$\Delta^{k_3} P^m(z)=\Delta^{k_3} P^{k_3+(k_1+k_2)}(z)$
is non-zero only if  
\begin{align}\label{MainResult-2-pe3}
k_3\leq N.
\end{align}

From the observations above and 
Eqs.\,(\ref{MainResult-2-pe1}), (\ref{MainResult-2-pe2}),
(\ref{MainResult-2-pe3})
it is easy to see that, $\Delta^m (f(z) P(z)^m)\neq 0$ 
only if $m=k_1+k_2+k_3 \leq d+N$. In other words,
$\Delta^m (f(z) P(z)^m)=0$ for any $m>d+N$. 
Hence Theorem \ref{MainResult-2} holds.
\epfv

Note that all results used in the proof above for 
the $(\Leftarrow)$ part of the theorem also hold 
for all HN formal power series. Therefore 
we have the following corollary.

\begin{corol}
Let $P(z)$ be a HN formal power series such that 
the vanishing conjecture holds for $P(z)$. 
Then, for any polynomial $f(z)$, 
we have $\Delta^m (f(z) P(z)^m)=0$ 
when $m>>0$.
\end{corol}

{\small \sc ${}^{*}$ Department of Mathematics, Radboud University Nijmegen, Postbus 9010,   
 6500 GL Nijmegen, The Netherlands.}

{\em E-mail}: essen@math.ru.nl

{\small \sc ${}^{**}$ Department of Mathematics, Illinois State University,
Normal, IL 61790-4520.}

{\em E-mail}: wzhao@ilstu.edu.

\end{document}